\documentclass[11pt]{amsart}
\usepackage{amssymb}
\usepackage{amsmath}
\usepackage{amscd}
\usepackage{amsthm}
\usepackage{verbatim}  
\usepackage[all]{xy}

\newtheorem{thm}	  {Theorem}      [section]
\newtheorem{theorem}   [thm] {Theorem}
\newtheorem{prop}   [thm] {Proposition}
\newtheorem{lemma}  [thm] {Lemma}

\newtheorem{corollary} [thm] {Corollary}
\theoremstyle{remark}
\newtheorem{remark} [thm] {Remark}
\newtheorem{definition}[thm]{Definition}

\DeclareMathOperator\aut{Aut}

\DeclareMathOperator\PSL{PSL}
\DeclareMathOperator\PSU{PSU}
\DeclareMathOperator\PGL{PGL}

\DeclareMathOperator\GL{GL}

\DeclareMathOperator\out{Out}
\renewcommand{\bar}[1]{#1\llap{$\overline{\phantom{\rm#1}}$}}

\newcommand\bG{\bar G}

\newcommand\cB{\mathcal{B}}
\newcommand\cC{\mathcal{C}}

\let\lift=\widehat

\newcommand\directsum{\oplus}
\newcommand\st{|}
\newcommand\tG{\tilde G}
\newcommand\tC{\tilde C}
\newcommand\tQ{\tilde Q}
\newcommand\normal{\triangleleft}

\DeclareMathOperator\Mat{M}

\DeclareMathOperator\PGU{PGU}
\newcommand{\inverselimit}{\mathop{\varprojlim}\limits}
\newtheorem{example}[thm]{Example}

\newcounter{parts1}
\newenvironment{parts1}
{
\begin{list}{{\bf
Part~\Alph{parts1}}}{\usecounter{parts1}\leftmargin-1pt} 
}
{\end{list}}

 \font\Aaa=msam10
        \def\nor{\hbox{\Aaa C}}

\begin{document}

\title
[Frobenius subgroups of free profinite products]
{Frobenius subgroups of free profinite products}

\author{Robert M. Guralnick}
\address{Department of Mathematics, University of Southern California,
Los Angeles, CA 90089--2532, USA}
\email{guralnic@usc.edu}

\author{Dan Haran}
\address{School of Mathematical Sciences,
Tel Aviv University,
TEL AVIV 69978, ISRAEL}
\email{haran@post.tau.ac.il}

\thanks{The first author was partially supported by NSF grant
DMS 0653873. 
The second author was partially supported by the
Israel Science Foundation (Grant No.~343/07)
and the Minerva Minkowski Center for Geometry at Tel Aviv University.
}

\begin{abstract}
We solve an open problem of Herfort and Ribes:
Profinite Frobenius groups of certain type
do occur as closed subgroups of
free profinite products of two profinite groups.
This also solves a question of Pop
about prosolvable subgroups of free profinite products.
\end{abstract}

\maketitle

\section{Introduction}

Herfort and Ribes show in
\cite[Theorem~3.2]{HR1}
that a closed solvable subgroup of the free product
of a family of profinite groups $\{A_x\}_{x \in X}$
must be one of the following:
\begin{enumerate}
\item
a conjugate of a subgroup of one of the free factors
$A_x$;
\item
isomorphic to 
$\hat{\mathbb Z}_\sigma \rtimes \hat{\mathbb Z}_{\sigma'}$,
where $\sigma$ and $\sigma'$ are disjoint sets of prime numbers;
\item
free pro-$\cC$ product of two copies of the group of order 2,
for some full class of $\cC$ of finite groups;
\item 
a profinite Frobenius group of the form
$\hat{\mathbb Z}_\sigma \rtimes C$
with Frobenius kernel 
$\hat{\mathbb Z}_\sigma$,
where $C$ is a finite cyclic group.
\end{enumerate}

In
\cite[Section~4]{HR1}
they show that
each group of the one of the first three types
does occur as a closed subgroup of a free profinite product,
namely, of two finite groups.
As for the Frobenius groups,
Herfort and Ribes state in \cite{HR1}
and show in \cite{HR2}
that they occur as closed subgroup of free {\it prosolvable} products
of two finite groups.
They (implicitly) leave open the question
whether the above Frobenius groups
occur as closed subgroups of free profinite products
of, say, two finite groups.

This problem has been explicitly posed in \cite{RZ}
(see \cite[Open Question~9.5.5]{RZ}).

The main result of this paper is 
an affirmative answer to this question.

The proof uses the classification of finite simple groups
(by analyzing automorphisms of simple groups and subgroups stabilized by them).
If $C$ has prime power order, the proof is much simpler
and does not require the classification
(Sylow's theorem is the main tool in that case).
This had been essentially done by the authors several years ago
in an unpublished work.
See Remark~\ref{prime power case}.
A.~Zalesski and P.~Zalesskii also found an independent proof in this case.

We thank P.~Zalesskii for reminding us
that our result also answers a question of Pop \cite{Pop}.
That paper characterizes
closed prosolvable groups of 
free profinite products of profinite groups
and the question is the following.
Is there a free profinite product 
$G = \coprod_{i \in I} G_i$
and a closed prosolvable subgroup $H$ of $G$
such that
\begin{enumerate}
\item
there is no prime $\ell$ such that
$H \cap G_i^\sigma$ is a pro-$\ell$ group
for all $i \in I$ and all $\sigma \in G$;
\item
$H \le G_i^\sigma$ for no $i \in I$ and no $\sigma \in G$?
\end{enumerate}
A Frobenius subgroup of a free product of certain
two finite groups provides such an example
(Example~\ref{Pop's question}).

\section{Finite and profinite Frobenius Subgroups}

Recall that the notions of order and index extend
from finite group to profinite groups;
instead of natural numbers these are supernatural numbers
(\cite[Section~2.3]{RZ} or \cite[Section~22.8]{FJ}).
In particular,
a profinite group has Sylow $p$-subgroups
for each prime $p$
(\cite[Corollary~2.3.6]{RZ} or \cite[Section~22.9]{FJ}).

A profinite group $F$ is a {\em Frobenius group}
if it is a semidirect product $F = C \ltimes K$
of nontrivial profinite groups $C,K$,
where $C$ acts on $K$ so that
$[c,k] \ne 1$ for every $1 \ne c \in C, 1 \ne k \in K$.
One then calls $K$ the {\em Frobenius kernel}
and $C$ a {\em Frobenius complement} of $F$.

Since we deal only with a special type of Frobenius groups,
we adopt the following notation.
Let $C$ be a finite cyclic group.
A {\em $C$-group}
is a profinite group with a distinguished subgroup
isomorphic to $C$; we identify this subgroup with $C$.
A {\em $C$-homomorphism} of $C$-groups $G \to H$
is a continuous homomorphism $G \to H$
that maps the copy of $C$ in $G$ identically onto the copy of $C$ in $H$.
If $A \nor G$ such that $A \cap C = 1$,
then $G/A$ is a $C$-group
(we identify $C$ with $CA/A$)
and $G \to G/A$ is a $C$-epimorphism.

We call a profinite $C$-group $F$ 
a {\em $C$-Frobenius group}
if $F$ is a Frobenius group with complement $C$ and procyclic kernel.
The following properties are easy to verify:

\begin{lemma} \label{properties of frobenius groups}
Let $C \ne 1$ be a finite cyclic group acting on 
a procyclic group $K \ne 1$.
Then
$F = C \ltimes K$ is a $C$-Frobenius group
if and only if
for each prime $p$ dividing the order of $K$,
the order of $C$ divides $p-1$
and $C$ acts faithfully on the $p$-primary part of $K$.
If $F$ is a $C$-Frobenius group then:
\begin{enumerate}
\renewcommand{\theenumi}{\alph{enumi}}
\item 
Every prime divisor of $|C|$ 
is strictly smaller than
any prime divisor of $|K|$.
\item 
$K$ is of odd order.
\item 
Any quotient group ${\bar F}$ of $F$ is
either a quotient of $C$
or a ${\bar C}$-Frobenius group, where
${\bar C}$ is the image of $C$ in ${\bar F}$.
\item  
Let ${\hat K} \to K$ be an epimorphism of procyclic groups
of orders divisible by the same primes.
Suppose that $C$ acts on ${\hat K}$
such that ${\hat K} \to K$ is $C$-equivariant.
Then $C \ltimes {\hat K}$ is also a $C$-Frobenius group.
\item  \label{special form of F_1} 
Every subgroup of $F$
is a conjugate of $C_1K_1$,
where $C_1 \le C$ and $K_1 \le K$.
\item  \label{normal subgroup} 
A subgroup of $F$ is normal if and only if
it is either a subgroup of $K$
or of the form $C_1K$, 
where $C_1 \le C$.
In particular, a minimal normal subgroup of $F$
is a minimal normal subgroup of $K$.

\item \label{g} 
Let $C_1 \le C$, $K_1 \le K$,
and $f \in F$.
Then 
$$
C_1^fK_1 \cap K = K_1
\qquad
\mbox{ and }
\qquad
C_1^fK_1 \cap C = 
\begin{cases}
C_1 & \text{if $f \in CK_1$} \cr
1   & \text{if $f \notin CK_1$.} \cr
\end{cases}
$$
\end{enumerate}
\end{lemma}

\begin{lemma} \label{orbit of a frobenius group}
Let $F = CK$ be a finite $C$-Frobenius group (with Frobenius kernel $K$).
Suppose $F$ acts transitively on a set $\Delta$.
Then
\begin{enumerate}
\item 
There is $L \in \Delta$
such that
its $F$-stabilizer is $C_1K_1$ with $C_1 \le C$ and $K_1 \le K$.
Fix such $L$. Then
\item 
$C_1$ is the $C$-stabilizer of $L$,
and hence also of every $L^c$, with $c \in C$.
\item
Every point of $\Delta \smallsetminus \{L^c \st c\in C\}$
has a trivial $C$-stabilizer.
\end{enumerate}
\end{lemma}

\begin{proof}
(a)
Let $L \in \Delta$.
Its $F$-stabilizer $F_1$ is,
by Lemma~\ref{properties of frobenius groups}(\ref{special form of F_1}),
a conjugate of $C_1 K_1$
for some $C_1 \le C$ and $K_1 \le K$.
Replacing $L$ by a conjugate we may assume
that $F_1 = C_1K_1$.

(b),(c)
By Lemma~\ref{properties of frobenius groups}(\ref{g})
$$
(L^f)^c = L^f
\iff
c \in (C_1K_1)^f
\iff
c \in C_1^fK_1 \cap C = 
\begin{cases}
C_1 & \text{if $f \in CK_1$}\cr
1   & \text{if $f \notin CK_1$}\cr
\end{cases}
$$
and $L^{CK_1} = L^{K_1C} = L^C$.
\end{proof}

\section{Intravariant Subgroups}

\begin{definition}
Let $H \le G$ be groups and let $A$ be a group acting on $G$
from the right.
We say that $H$ is {\em $A$-intravariant in $G$}
if for every $x \in A$ 
there is $g \in G$ such that $H^x = H^g$.
We say that $H$ is {\em an intravariant subgroup of $G$}
if it is $\aut(G)$-intravariant in $G$.
\end{definition}

We point out that Sylow subgroups and their normalizers
are intravariant subgroups.
In the rest of this section we exhibit further families
of intravariant subgroups of finite simple groups.

Recall that an {\em almost simple group} is a group $G$
with a unique minimal normal subgroup $S$
which is a nonabelian simple group.
Thus, $S \nor G \le \aut(S)$.
We refer
the reader to \cite{AS, C, GL, GLS3} for the basic facts
about automorphisms of finite simple groups. 

We recall some facts about automorphisms of the finite
simple groups and most especially about Chevalley groups.
The most complicated cases to deal with are $\PSL$ and $\PSU$.

Our first result is \cite[3.22]{GMS}.

\begin{lemma} \label{centralizer of an automorphism is not trivial}
Let $S$ be a finite nonsolvable group.
Let $x \in \aut(S)$.
Then $C_S(x) \ne 1$. 
\end{lemma}


The next result is the Borel-Tits Theorem \cite[Theorem 3.1.3]{GLS3}.

\begin{lemma} \label{borel-tits}
Let $S$ be a simple Chevalley group
and $U$ a nontrivial unipotent subgroup.
Let $A=\aut(S)$.
Then there exists a proper parabolic
subgroup $P$ of $S$ such that $U$ is contained in the unipotent radical
of $P$ and $N_A(P) \ge N_A(U)$.
\end{lemma}

We remark that $P$ is proper, since the unipotent radical of $P$
is normal in $P$, while $S$ is simple.

We also require:

\begin{lemma} \label{sylow}
Let $S \nor G$ be finite groups.
Let $R$ be an
$r$-subgroup of $G$ for some prime $r$.
\begin{enumerate}
 \item \label{r divides S}
If $r$ divides $|S|$, then $R$ normalizes some Sylow $r$-subgroup of $S$.
\item \label{r does not divide S}
If $r$ does not divide $|S|$ but another prime $p$ does divide $|S|$,
then $R$ normalizes a Sylow $p$-subgroup of $S$.
\end{enumerate}
\end{lemma}

\begin{proof}
Let $Q$ be a Sylow $r$-subgroup of $G$ containing $R$.
Then $R$ normalizes $Q \cap S$,
which gives (\ref{r divides S}).

To prove (\ref{r does not divide S}),
let $P$ be a Sylow $p$-subgroup of $S$.
By Sylow's theorem,
$G=SN_G(P)$.
Since $r$ does not divide $|S|$, $N_G(P)$ 
contains a Sylow $r$-subgroup of $G$,
say, $Q^g$, with $g \in G$.
Then $R \le Q \le N_G(P^{g^{-1}})$,
that is, $R$ normalizes $P^{g^{-1}}$.
\end{proof}

We now examine $\PSL$ and $\PSU$ more closely.
Let $S=\PSL(d,p^e)$, where $p$ is a prime and $S$ is simple.
Let $\sigma$ be the Frobenius automorphism
and $\tau$ the graph automorphism
(which we may view as the inverse transpose map) of $S$.
Then
$$
\Omega = \PGL(d,p^e) \cup\{\sigma,\tau\}
$$
generates $A = \aut(S)$.
Moreover,
$\PGL(d,p^e)\langle \sigma \rangle$ is of index $2$ in $A$.

We will use the following elementary result from linear algebra.

\begin{lemma}\label{determinant is norm}
Let $F$ be a field.
Let $A \in \Mat_n(F)$ such that its minimal polynomial $g$ is irreducible
of degree $d$.
Then $K = F[A]$ is a finite field extension of $F$
and for every $A' \in K$ we have
$\det A' = N^K_F(A')^{n/d}$,
where
$N^K_F \colon K^\times \to F^\times$ is the norm.
\end{lemma}

\begin{proof}
As $g$ is irreducible,
$K \cong F[X]/(g)$ is a field extension of $F$ of degree $d$.

For each $m \in \mathbb {N}$ and $A' \in K$
let $A'_m \colon K^m \to K^m$ be the $F$-linear map
given by $A'_m(x_1,\ldots, x_m) = (A'x_1,\ldots, A'x_m)$.
By definition,
$N^K_F(A')$ is the determinant of $A'_1 \colon K \to K$.
Hence
$N^K_F(A')^m = \det A'_m$.

Now,
the matrix of $A_1$ with respect to the basis
$I,A,\ldots, A^{d-1}$ of $K$ is the companion matrix $C$ of $g$.
Thus the matrix of $A_m$ with respect to some basis of $K^m$ 
is the direct sum $C^{(m)}$ of $m$ copies of $C$.
But $A$ is similar $C^{(m)}$,
where $m = \frac{n}{d}$
(this is the rational canonical form of $A$),
hence
$A$ is the matrix of $A_m$ with respect to a suitable basis $\cB$ of $K^m$.

As $A' \in K = F[A]$,
there is $h \in F[X]$ such that $A' = h(A)$.
Clearly $A'_m = h(A)_m = h(A_m)$ and
the matrix of $h(A_m)$ with respect to $\cB$ is $h(A) = A'$.
Hence
$N^K_F(A')^m = \det A'_m = \det (A')$.
\end{proof}

\begin{lemma}\label{det of centralizer of semisimple}
Let $F$ be a finite field.
Let $B \in \Mat_d(F)$ be semisimple.
Let $a \in F$.
Then there is $B' \in \Mat_n(F)$ 
such that $B B' = B' B$ and $\det(B') = a$.
\end{lemma}

\begin{proof}
Without loss of generality $B$ is in rational canonical form.
Thus $B = {\rm diag}(B_1,\ldots, B_r)$,
where each $B_i \in \Mat_{n_i}(F)$
has minimal polynomial irreducible of degree $n_i$.
Then by Lemma~\ref{determinant is norm}
$K = F[A]$ is a finite field extension of $F$.
Hence the norm $N^K_F \colon K^\times \to F^\times$
is surjective.
Therefore, again by Lemma~\ref{determinant is norm},
there is $B'_i \in \Mat_{n_i}(F)$ such that
$B'_i \in F[B_i]$, whence $B'_i$ commutes with $B_i$,
and $\det(B'_1) = a$ and $\det(B'_i) = 1$ for $i >1$.
Then $B = {\rm diag}(B'_1,\ldots, B'_r)$
commutes with $B$ and $\det(B') = a$.
\end{proof}

\begin{lemma} \label{intravariant subgroup}
Let $S=\PSL(d,p^e)$ and $H=\PGL(d,p^e)$.
Put $A = \aut(S)$.
\begin{enumerate}
\item\label{parabolic is intravariant}
Any parabolic subgroup $P$ of $S$
whose normalizer contains 
an element outside of $H\langle \sigma \rangle$ is intravariant.
\item \label{semisimple is intravariant}
If $h \in H$ is a semisimple element
(i.e. has order prime to $p$),
then $\langle h \rangle$ is $A$-intravariant in $H$.
\item \label{normalizer and centralizer of semisimple is intravariant}
If $h \in H$ is semisimple,
then $C_S(h)$ and $N_S(\langle h \rangle)$
are intravariant subgroups of $S$.
\end{enumerate}
\end{lemma}

\begin{proof}
(\ref{parabolic is intravariant})
We may conjugate $P$ in $S$ and thus assume that
$P$ contains the standard Borel subgroup of $S$
that consists of upper triangular matrices
of determinant $1$.
Clearly, diagonal matrices
and $\sigma$ normalize $P$.
Hence $P$ is an intravariant subgroup of $S$ if and only if
$N_A(P)$ contains
an element outside of $H\langle \sigma \rangle$.

(\ref{semisimple is intravariant})
Let $x \in A$.
We have to show that 
there is $s \in H$ such that
$\langle h \rangle^x = \langle h \rangle^s$.
We may assume that $x \in \Omega$.
If $x \in H$, the assertion is trivial.
So assume that
either $x = \sigma$ or $x = \tau$.
Put $F = \mathbb{F}_{p^e}$.
Lift $h$ to a semisimple element of $\GL(d,F)$.
It suffices to show that
there is $m \in \mathbb{Z}$ such that
$h^{x}, h^m$ are conjugate in $\GL(d,F)$,
i.e., similar over $F$.

Since every square matrix is similar over $F$ to its transpose,
$h^\tau$ is similar to $h^{-1}$.
On the other hand,
$h^\sigma$ is similar over $F$ to $h^p$.
Indeed,
consider $\GL(d,F)$ as a subgroup of $\GL(d, \bar F)$
and 
extend $\sigma$ to the Frobenius automorphism of $\GL(d,\bar F)$.
There is $z \in \GL(d, \bar F)$
such that $h^z \in \GL(d, \bar F)$ is diagonal.
Then clearly $(h^z)^\sigma = (h^z)^p$.
Therefore
$(h^\sigma)^{z^\sigma} = (h^z)^\sigma = (h^z)^p = (h^p)^z$.
Hence $h^\sigma, h^p$ are similar over $\bar F$.
Therefore they are similar over $F$.

(\ref{normalizer and centralizer of semisimple is intravariant})
Let $x \in A$.
By (\ref{semisimple is intravariant}) there is $s \in H$ such that
$\langle h \rangle^x = \langle h \rangle^s$.
By Lemma~\ref{det of centralizer of semisimple}
there is $z \in H$ such that $zh = hz$ and $\det(z) \equiv \det(s)
\pmod{(F^\times)^d}$.
Replace $s$ by $z^{-1}s$ to get that $s \in S$.
As $C_S(h) = S \cap C_H(h)$,
$N_S(h) = S \cap N_H(h)$,
and $S^x = S = S^s$,
we get 
$C_S(h)^x = C_S(h)^s$
and
$N_S(h)^x = N_S(h)^s$.
%
\end{proof}

\begin{prop} \label{PSL}
Let $S=\PSL(n,p^e)$ be simple.
Then every $x \in \aut(S)$
normalizes a nontrivial proper intravariant subgroup of $S$.
\end{prop}

\begin{proof}
Let $H=\PGL(n,p^e)$.
If $C_H(x)$ is not a $p$-group,
we may choose $1 \ne h \in C_H(x)$ of order prime to $p$.
Thus, $h$ is semisimple and $C_S(h)^x = C_S(h^x) = C_S(h)$.
By Lemma~\ref{intravariant subgroup}%
(\ref{normalizer and centralizer of semisimple is intravariant}),
$C_S(h)$ is an intravariant subgroup of $S$.
By Lemma~\ref{centralizer of an automorphism is not trivial},
$C_S(h) \ne 1$;
clearly $C_S(h) \ne S$.
So we may assume
that $C_H(x)$ is a $p$-group (i.e. consists of unipotent elements).
In particular, $C_H(x) \le S$, whence $C_S(x) = C_H(x)$.

If $x \in \aut(S) \smallsetminus H\langle \sigma \rangle$,
then, by Lemma~\ref{borel-tits},
$x$ normalizes some proper parabolic subgroup of $S$
which is intravariant by Lemma~\ref{intravariant subgroup}%
(\ref{parabolic is intravariant}).

If $x$ is in $H$,
then $x \in C_H(x) \le S$,
so $x$ is contained in some Sylow $p$-subgroup of $S$,
which is intravariant.

So we may assume that $x \in H\langle \sigma \rangle \smallsetminus H$.
Write $x = yz = zy$,
where $y \in \langle x \rangle$ has order prime to $p$
and $z \in \langle x \rangle$ has order a power of $p$.
The restriction to $\langle y \rangle$
of the projection $H\langle \sigma \rangle \to \langle \sigma \rangle$
is injective,
because
its kernel $H \cap \langle y \rangle$
is contained in the $p$-group $C_H(x)$.
Thus, writing $y$ as $\sigma^j h$, with some $h \in H$,
we see that $y$ and $\sigma^j$ have the same order.
By Shintani descent
\cite[p.~81]{GL}
they are conjugate by an element of $H$.
Since $C_H(\sigma^j)$ contains the Jordan matrix $J_n(1)$,
it follows that $C_H(y)$ contains a regular unipotent element
(a conjugate of $J_n(1)$).

As $z$ commutes with $y$, it normalizes $C_H(y)$.
Its order is a power of $p$,
and so $z$ normalizes some Sylow $p$-subgroup $T$ of $C_H(y)$.
But $y$ centralizes $T$
and so $x = yz$ normalizes $T$.
Since $T$ contains a regular unipotent element,
$T$ is contained in a unique Sylow $p$-subgroup $P$ of $H$.
As $T = T^x \le P^x$, we have $P^x = P$.
Thus $x$ normalizes this intravariant subgroup.
\end{proof}

The same result holds for $S=\PSU(d,p^e)$.

\begin{lemma} \label{PSU}
Let $S=\PSU(d,p^e)$ be simple.
If $x \in \aut(S)$, then
$x$ normalizes a nontrivial proper intravariant subgroup of $S$.
\end{lemma}

\begin{proof}
In this case
$\aut(S) = H \langle \sigma\rangle$,
where
$H = \PGU(d,p^e)$
and $\sigma$ is the Frobenius automorphism of order $2e$.
Let $x \in \aut(S)$.
If $C_H(x)$ contains a nontrivial semisimple element $h$, then
$C_S(h)$ is an intravariant subgroup normalized by $x$.
Otherwise, $C_H(x)$ is unipotent and by Lemma~\ref{borel-tits},
$x$ normalizes some (proper) parabolic subgroup of $S$.
In this case, all parabolics are intravariant, whence the result.
\end{proof}

The analogous result is true
for almost all of the families of simple groups.
However, it is not true for $S=\Omega^+(8,q)$.
Take $x \in S$ of order $(q^4-1)/\gcd(2,q-1)$.
It is not difficult to see that $x$ is contained
in no maximal subgroup of $\aut(S)$ not containing $S$,
whence it cannot normalize any nontrivial intravariant subgroup of $S$.
If we only consider cyclic groups of automorphisms of simple groups
which do not contain any inner automorphisms, then the result is true.
We will revisit this topic in a future paper.
It is not required for the results of this paper.

The outer automorphism groups of the other simple groups
are less complicated.
It is convenient to use the notation
$\PSL^{\epsilon}(d,q)$ with $\epsilon = \pm$---here
$\PSL^{+}$ means  $\PSL$ and $\PSL^{-}$ means $\PSU$.

\begin{lemma} \label{frobenius quotient}
Let $G$ be an almost simple finite group with socle $S$.
Suppose that $G/S$ is a Frobenius group with
a cyclic Frobenius kernel $K$ and
a cyclic Frobenius complement $C$.
If $|C| > 3$, then
$S \cong \PSL^{\epsilon}(d,q)$
and the Frobenius kernel is contained
in the subgroup of diagonal automorphisms of $\out(S)$.
Moreover, $d \ge 5$ and $q > 5$.
\end{lemma}

\begin{proof}
We have $G/S \le \out(S)$.
If $S$ is an alternating or sporadic group,
then $\out(S)$ and hence also $G/S$ has exponent $2$,
a contradiction to $G/S$ being Frobenius.
So $S$ is a Chevalley group. 

Let $r$ be a prime dividing $|K|$.
By Lemma~\ref{properties of frobenius groups},
$|C|$ divides $r-1$.
As $|C| > 3$, we have $r \ge 5$.
As $K$ is cyclic,
$C$ normalizes its Sylow $r$-subgroup.
If $S$ is not $\PSL^{\epsilon}(d,q)$,
then the Sylow $r$-subgroup of $K$ consists of field automorphisms.
However, the normalizer in $\out(S)$ of a field automorphism
is its centralizer.
Thus, $C$ centralizes a nontrivial subgroup of $K$,
a contradiction to $G/S = C\ltimes K$ being Frobenius.
This proves that $S\cong \PSL^{\epsilon}(d,q)$.
Arguing similarly,
it follows that $K$ consists of diagonal automorphisms and so $r|\gcd(d,q - \epsilon 1)$,
whence $d \ge 5$ and $ q > 5$.
\end{proof}

We can now show:

\begin{corollary}
\label{simple-reduct}
Let $G$ be an almost simple finite group with socle $S$.
Suppose that $G/S$ is a Frobenius group
with cyclic Frobenius kernel $K$ and cyclic Frobenius complement $C$.
Let $D$ be any cyclic subgroup of $G$ with $DS/S=C$.
If $C$ has prime power order, assume the same is true of $D$.
Then $D$ normalizes
a proper nontrivial intravariant subgroup $H$ of $S$.
\end{corollary}

\begin{proof}
If $|C|$, and hence also $|D|$, is a prime power,
then $D$ normalizes a nontrivial Sylow subgroup of $S$
by Lemma~\ref{sylow}.
Otherwise, $|C| > 3$,
and hence by Lemma~\ref{frobenius quotient},
$S$ is either $\PSL$ or $\PSU$.
Now apply Proposition~\ref{PSL} and Lemma~\ref{PSU}.
\end{proof}

\section{Lifting Frobenius Groups}

To show our main result we need some preparations.

\begin{lemma} \label{fiber product}
Let $\rho_1 \colon F_1 \to F_3$
and $\rho_2 \colon F_2 \to F_3$
be $C$-epimorphisms of $C$-Frobenius groups.
Then $F_1 \times_{F_3} F_2$ contains a $C$-Frobenius group
mapped onto $F_1$.
\end{lemma}

\begin{proof}
Let $K_i$ be the Frobenius kernel of $F_i$, 
and $K_i^{(p)}$
its Sylow $p$-subgroup, respectively,
for $i = 1,2,3$ and for each prime $p$.
Put $F = F_1 \times_{F_3} F_2$.
We identify $C$ with $C \times_{F_3} C$
and thus $F$ is a $C$-group
and the coordinate projections
$F \to F_1$, $F\to F_2$ are $C$-epimorphisms.
Clearly every $1 \ne c \in C$ acts
fixed-point-freely on the abelian subgroup
$K = K_1 \times_{F_3} K_2 $ of $F$.
For each prime $p$
we find below a cyclic $p$-subgroup
$\langle k_p \rangle$ of $K$
normalized by $C$
and mapped by $F \to F_1$ onto $K_1^{(p)}$.
Then
$C \big( \prod_p \langle k_p \rangle \big)$
is a $C$-Frobenius group mapped onto $F_1$.

Fix  a generator $k_1$ of $K_1^{(p)}$
and put $k_3 = \rho_1(k_1)$.
By Sylow's theorem,
$\rho_1(K_1^{(p)}) = K_3^{(p)} = \rho_2(K_2^{(p)})$.
So there is $k_2 \in K_2^{(p)}$
such that $k_3 = \rho_2(k_2)$.
If $k_3 = 1$, take $k_2 = 1$;
otherwise $k_2$ generates $K_2^{(p)}$.
Then $k = (k_1,k_2) \in F$
and $\langle k \rangle \le F$ is normalized by $C$,
that is,
$(k_1^c,k_2^c) \in \langle k \rangle$,
where $c$ is a generator of $C$.

Indeed,
if $k_2 = 1$, the assertion is clear.
So assume that $k_2 \ne 1$ and hence $k_3 \ne 1$.

Let $i = 1,2, 3$.
Writing 
$K_i^{(p)}$
additively,
it is a quotient of $\mathbb{Z}_p$
so that
$k_i$ is the class of the generator $1$ of $\mathbb{Z}_p$
and
$\rho_1 \colon K_1^{(p)} \to K_3^{(p)}$, 
$\rho_2 \colon K_2^{(p)} \to K_3^{(p)}$
are the quotient maps.
As the order of $c$ is prime to $p$,
there is a unique $m_i \in (\mathbb{Z}_p)^\times$
such that $c$ acts on $K_i^{(p)}$ by multiplication by the image of $m_i$.
But $\rho_1, \rho_2$ are $C$-equivariant,
hence $m_1 = m_3 = m_2$.
Thus
$(k_1^c,k_2^c) = (m_1 k_1, m_2 k_2) = m_1 k \in \langle k \rangle$,
\end{proof}

\begin{lemma}
\label{A is elementary abelian}
Let $G$ be a finite group and $A$ a minimal normal subgroup of $G$.
Assume that
$A$ is an elementary abelian $p$-group
and $G/A$ is a $C$-Frobenius group.
Then $G$ contains a $C$-Frobenius group $H$ such that $G = HA$.
\end{lemma}

\begin{proof}
Put $F = G/A$ and let $K$ be the Frobenius kernel of $F$.
Let $A \le M \normal G$ such that $K = M/A$.
So $G = C \ltimes M$.
For each prime $r$ let $K_r$ be the
Sylow $r$-subgroup of $K$
and let $M_r$ be a Sylow $r$-subgroup of $M$
such that $M_rA/A = K_r$.
As $A$ acts trivially on itself, $F$ acts on $A$,
and $A$ is an irreducible $F$-module.

Assume first that $K$ acts nontrivially on $A$.
Then there is a prime $r$ such that $K_r$ acts nontrivially on $A$.
Thus $A^{K_r} \ne A$.
As $K_r$ is normal in $F$,
$A^{K_r}$ is an $F$-submodule of $A$,
and hence $A^{K_r} = 0$.
Therefore $K_r$ 
has no nontrivial fixed points on $A$.
In particular, $r \ne p$ and $K_r \ne 1$.
Thus $M_r$ maps modulo $A$ isomorphically onto $K_r$
and $C_A(M_r) = 1$.
As $M_r \cap A = 1$,
we have $N_A(M_r) = C_A(M_r)$.
%
So $N_A(M_r) = 1$.
By Sylow's theorem, $G=AN_G(M_r)$.
%
%
Thus, $N_G(M_r)$ is a complement to $A$ in $G$,
and hence isomorphic to $F$.
In particular $N_G(M_r)$ is a $C_1$-Frobenius group for some
cyclic subgroup $C_1$,
isomorphic to $C$.
Note that $A$ (after scalars extension)
is the direct sum of eigenspaces $A_\lambda$ of $K_r$
corresponding to eigenvalues $\lambda$ of $r$-power order.
Since generators $k$ of $K_r$ and $c$ of $CA/A$ satisfy
$kc = ck^n$, 
where $n = |C|$ is prime to $r$,
it follows that $C$ maps $A_\lambda$ onto $A_{\lambda^n}$.
Therefore $C$ permutes freely the $A_\lambda$'s.
Thus, $A$ is a free $C$-module.
In particular, $H^1(C,A)=0$ and so
any two complements to $A$ in $AC$ are conjugate.
Thus, $C$ is conjugate to $C_1$
and so is contained in a complement $H$ to $A$.

Next assume that $K$ acts trivially on $A$.
As $K$ is cyclic, $M$ is abelian.
It suffices to find, for each prime $r$ dividing $|K|$,
a cyclic $r$-subgroup $L_r$ of $M$,
normalized by $C$,
such that $L_rA/A = K_r$.
Indeed, then
$H = (\prod_r L_r) C$
is a $C$-Frobenius group mapped onto $F$.

Since $M$ is abelian,
$M_r$ is normalized by $C$.
So if $M_r$ is cyclic, take $L_r = M_r$.
This is the case if $r \ne p$,
since $M_r \cong K_r$.
So let $r=p$ with $M_p$ not cyclic and $K_p \ne 1$.
Thus, $|C|$ divides $p-1$.
Hence every vector space over $\mathbb{Z}/p\mathbb{Z}$ 
on which $C$ acts is the direct sum of $1$-dimensional $C$-modules.
Since $C$ acts irreducibly on $A$,
this implies that $A$ is $1$-dimensional.  
As $K_p$ is cyclic and $M_p$ not, $M_p$ is of rank $2$,
whence its Frattini quotient $V$ is of dimension $2$.
In particular, the image ${\bar A}$ of $A$ in $V$ is of dimension $1$.
Thus $V = {\bar A} \directsum {\bar B}$
for some $1$-dimensional $C$-module ${\bar B}$.
Choose a cyclic subgroup $L_p$ of $M_p$ such that
its image in $V$ is ${\bar B}$.
Then, $L_p$ is normalized by $C$
and $M_p = L_p A$, whence $L_pA/A = K_p$.
\end{proof}

\begin{lemma}
\label{intravariant subgroup of a product}
Let $G$ be a finite group and let $A \normal G$.
Assume that
$A = \prod_{i=1}^t Q_i$
and the conjugation in $G$ transitively permutes the $Q_i$.
Let $G_1 = N_G(Q_1)$
and for each right coset $Z$ of $G_1$ in $G$
let $\lift{Z} \in Z$.
Let $U_1$ be a $G_1$-intravariant subgroup in $Q_1$.
Then 
$U := \prod_Z U_1^{\lift{Z}} \le A$
is $G$-intravariant in $A$.
Moreover, if $U_1$ is not normal in $Q_1$
then $U$ is not normal in $A$.
\end{lemma}

\begin{proof}
We have
$\{Q_i\}_{i=1}^t = \{Q_1^{\lift{Z}}\}_{Z \in G/G_1}$.
Thus $A = \prod_{Z \in G/G_1} Q_1^{\lift{Z}}$
and hence $U \le A$.

To show the intravariance,
let $g \in G$.
For each $Z \in G/G_1$ we have
$Q_1^{ \lift{Z}g\lift{Zg}^{-1} } = Q_1$,
hence $\lift{Z}g\lift{Zg}^{-1} \in G_1$.
Thus there is $b_Z \in Q_1$ such that
$U_1^{ \lift{Z}g\lift{Zg}^{-1} } = U_1^{b_Z}$.
Put $a_Z = \lift{Zg}^{-1} b_Z \lift{Zg}$,
then $a_Z \in Q_1^{\lift{Zg}}$
and
$U_1^{ \lift{Z}g } = (U_1^{\lift{Zg}})^{a_Z}$.
Thus
$a : = \prod_{Z \in G/G_1} a_Z \in A$
and $U^{g} = U^a$.

The last assertion is clear.
\end{proof}

\begin{lemma}
\label{invariant subgroup of a product}
Assume, in the situation of the preceding lemma,
that $G$ is a $C$-group such that $G/A$ is a $C$-Frobenius group.
Let $C_1 = N_C(Q_1) = C \cap G_1$
and assume that $G_1/A = C_1 K_1$
for some subgroup $K_1$ of the Frobenius kernel of $G/A$.
Assume that $U_1$ is $C_1$-invariant.
Then we may choose the $\lift{Z} \in Z$ so that
$U$ is $C$-invariant.
\end{lemma}

\begin{proof}
It suffices to choose the representatives $\lift{Z}$ so that
\begin{equation}
U_1^{\lift{Zc}} = U_1^{\lift{Z}c}
\hbox{ for every } c \in C .
\label{lifting}
\end{equation}
But to achieve~(\ref{lifting}),
it suffices to achieve it
only for a representative $Z$ of each $C$-orbit in $G/G_1$
(when $G$ and $C \le G$ act on $G/G_1$ by multiplication from the right).
Let $C_Z$ be the $C$-stabilizer of $Z$
and let $R_Z$ be a set of representatives of $C/C_Z$ in $C$.
Suppose we have found $\lift{Z} \in Z$
such that $C_Z$ normalizes $U_1^{\lift{Z}}$.
Then put
$\lift{Zcr} = \lift{Z}r$ for all $c \in C_Z$ and $r \in R_Z$
to get~(\ref{lifting}).

As $G$ acts transitively on $G/G_1$ 
with $A$ acting trivially,
this induces a transitive action of $F = G/A$ on $G/G_1$
and the $F$-stabilizer of $G_1$ is $C_1K_1$.
By Lemma~\ref{orbit of a frobenius group}
either $Z = G_1 c$ for some $c  \in C$ and $C_Z = C_1$
or $C_Z = 1$.
In the former case
we may take $\lift{Z} = c$
to make $C_Z$ normalize $U_1^{\lift{Z}}$.
In the latter case
we may choose $\lift{Z}$ arbitrarily.
\end{proof}

We now prove the main step:

\begin{theorem} \label{main finite}
Let $G$ be a finite group with a normal subgroup $A$
and a cyclic subgroup $C$
such that
$C \cap A=1$ and
$G/A$ is a $C$-Frobenius group.
Then $G$ contains
a $C$-Frobenius group $H$
such that $G = HA$.
\end{theorem}

\begin{proof}
Put $F = G/A$ and let $K$ be its Frobenius kernel.
We identify $CA/A$ with $C$.

We prove the result by induction on the order of $A$.
We may assume that $A \ne 1$.
We divide the proof into several parts.
\begin{parts1}
\item\label{kernel is minimal normal}
{\em We may assume that $A$ is a minimal normal subgroup of $G$.}
Let $B$ be a minimal normal subgroup $G$ contained in $A$.
By induction, $G/B$ satisfies the theorem with $A/B \normal G/B$.
Thus there is a subgroup $H_0$ of $G$ containing $B$ and $C$
such that
$H_0A=G$ and $H_0/B$ is a $C$-Frobenius group.
If $B \ne A$,
apply the induction to $H_0$ with $B \normal H_0$ 
to get
a $C$-Frobenius subgroup $H$ of $H_0$
such that $H_0 = HB$.
Then $HA = H_0A = G$.
So we may assume that $B = A$.

\item\label{A is a product of simple groups}
Thus $A$ is the direct product of copies of a finite simple group.
If $A$ is an elementary abelian $p$-group with $p$ prime,
then we are done by Lemma~\ref{A is elementary abelian}.
So we may assume that $A=Q_1 \times \ldots \times Q_t$ where $Q_i=Q$
is a nonabelian simple group.

\item\label{faithful action}
{\em We may assume that $G$ acts faithfully on $A$,
that is, $G \to \aut(A)$ is an embedding.}
Indeed, 
otherwise let $B$ be a minimal normal subgroup of $G$
contained in the kernel of this map.
Then $A \cap B = 1$,
and hence $B$ is isomorphic to a subgroup of $G/A$.
In particular, $B$ is solvable.
Moreover, the image of $B$ in $G/A$ is a minimal normal subgroup,
and hence is contained in $K$ by
Lemma~\ref{properties of frobenius groups}(\ref{normal subgroup}).
Thus, $AB \cap C=1$.
This allows us to identify $C$ with $CAB/AB$.

The quotient $G/AB$ of $G/A$
either equals to $C$
(that is, $G = ABC$)
or is a $C$-Frobenius group with $|G/AB| < |G/A|$.
In both cases
there is a solvable subgroup $G_0$ of $G$ containing $C$ and $B$,
such that $G_0A/B = G/B$:
In the first case take $G_0 = CB$.
In the second case
proceed by subinduction on $|G/A|$;
by the hypothesis,
there is a subgroup $G_0$ of $G$ containing $C$ and $B$,
such that $G_0A/B = G/B$ and $G_0/B$ is $CB/B$-Frobenius.
In particular, $G_0/B$ is solvable,
hence so is $G_0$.
In both cases,
as $G_0 A = G$,
we may replace $G$ by $G_0$ and $A$ by $A \cap G_0$.
As $G_0$ is solvable but $A$ not, $|A \cap G_0| < |A|$.
So the existence of $H$ follows by induction.

\item\label{reduction to intravariance}
{\em Reduction to intravariance}
Our aim is to construct a subgroup $U$ of $A$
such that
$U$ is $A$-intravariant in $G$,
$C$ normalizes $U$,
and $A$ does not normalize $U$.
Then
$N_G(U)A = G$
and $N_G(U) \cap A$ is a proper subgroup of~$A$.
Hence by induction hypothesis
$N_G(U)$ contains a $C$-Frobenius subgroup $H$
with $HA = G$.

\item\label{division}
{\em Division into three cases.}
So let $G_1 = N_G(Q_1)$.
Then $A \le G_1$; put $F_1 = G_1/A \le F$.
Note that
$F= G/A$ acts transitively on
the set $\Delta:=\{Q_1, \ldots, Q_t\}$,
and $F_1$ is the stabilizer of $Q_1$.
By Lemma~\ref{properties of frobenius groups}(\ref{special form of F_1})
we may replace $Q_1$ by some conjugate $Q_i$
to assume that $F_1 = C_1K_1$,
where $C_1 \le C$ and $K_1 \le K$.
We divide the rest of the proof into three cases.

{\em Case I: $C_1$ centralizes $Q_1$.}
In this case let $U_1 \ne 1$ be a Sylow subgroup of $Q_1$.
Then $U_1$ is $C_1$-invariant, intravariant in $Q_1$,
and not normal in $Q_1$.
By Lemma~\ref{intravariant subgroup of a product}
and Lemma~\ref{invariant subgroup of a product}
there is a  $C$-invariant $G$-intravariant subgroup $U$ of $A$
which is not normal in $A$.
(Notice that case $C_1 = 1$ is included here.)

{\em Case II: $C_1, K_1 \ne 1$.}
In this case $F_1$ is a $C_1$-Frobenius group.
Let $\phi \colon G_1 \to \aut(Q_1)$
be the map induced by the action of $G_1$ on $Q_1$
and let $\bG_1 = \phi(G_1)$.
Then $\phi(A) = \phi(Q_1) = Q_1$,
hence $\phi$ induces a surjection
$\bar\phi \colon F_1 = \bG_1/A \to \bG_1/Q_1 \le \out(Q_1)$.
We first claim that $\bar\phi$ is an isomorphism,
that is,
$F_1 \to \out(Q_1)$ is injective. 

It suffices to show that
$K_1 \to \out(Q_1)$ is injective,
since the kernel of $F_1 \to \out(Q_1)$ 
is a normal subgroup of the Frobenius group $F_1$,
which does not intersect its Frobenius kernel $K_1$
and hence is trivial
by Lemma~\ref{properties of frobenius groups}(\ref{normal subgroup}).

So let $g \in G$ with image in $K_1$
act on $Q_1$ as an inner automorphism.
As $F = G/A$ is a Frobenius group,
for every $\sigma \in G$ we have
$\sigma g \sigma^{-1} = g^m a$ for some $m \in \mathbb{N}$ and $a \in A$.
It follows that $g$ acts as an inner automorphism,
namely $(g^\sigma)^m a^\sigma$,
on $Q_1^\sigma$.
Thus $g$ acts as an inner automorphism on $A$.
As $G$ acts faithfully on $A$,
this means that $g \in A$.
Therefore the image of $g$ in $K_1$ is trivial.
This proves the claim.

Thus $\bG_1/Q_1 \cong F_1$ is a $C_1$-Frobenius group.
By Corollary \ref{simple-reduct}, 
$C_1$ normalizes some nontrivial intravariant proper
subgroup $U_1$ of $Q_1$.
So again we are done
by Lemma~\ref{intravariant subgroup of a product}
and Lemma~\ref{invariant subgroup of a product}.

{\em Case III: $K_1 = 1$ and $C_1$ does not centralize $Q_1$.}
In particular, $C$ does not centralize $Q_1$.
Notice that $G_1 = C_1 A \le CA$.
For each $g \in G$ let
$\Delta_g = \{Q_1^{\sigma g} \st \sigma \in CA \} \subseteq \Delta$.
Then
$\Delta_{g_1} \cap \Delta_{g_2} \ne \emptyset
\implies
CAg_1 = CAg_2
\implies
\Delta_{g_1} = \Delta_{g_2}$.
Thus $\Delta$ is the disjoint union of the distinct $\Delta_g$.
Therefore if we define
$\tQ_g = \langle Q_1^{\sigma g} \st \sigma \in CA \rangle$,
then $A$ is the direct product of the distinct $\tQ_g$.
The conjugation in $G$ permutes the $\Delta_g$
and hence also the $\tQ_g$.
Put $\tG_1 = N_G(\tQ_1)$; then $\tG_1 = CA$.
Thus $\tC_1 = C \cap \tG_1 = C$.

We define $U_1 = C_{\tQ_1}(\tC_1)$;
thus $U_1$ is a $\tC_1$-invariant subgroup of $\tQ_1$.
By our assumption, $U_1 \ne \tQ_1$.
By Lemma~\ref{centralizer of an automorphism is not trivial},
$U_1 \ne 1$. As
$\tG_1 = \tC_1 A =
\tC_1(\tQ_1 \times \prod_{\tQ_g \ne \tQ_1} \tQ_g)$
and $U_1$ is $\tC_1$-invariant,
$U_1$ is $\tG_1$-intravariant in $\tQ_1$.
Now apply again
Lemma~\ref{intravariant subgroup of a product}
and Lemma~\ref{invariant subgroup of a product}
with $\{\tQ_g \st g \in G\}$ instead of  $\Delta$
and with $\tQ_1,\tG_1,\tC_1$ instead of $Q_1, G_1, C_1$.  
\end{parts1}
\end{proof}

\begin{remark}\label{prime power case}
If $C$ is an $r$-group for some prime $r$,
we may considerably simplify the proof
of Theorem~\ref{main finite},
leaving out the classification of finite simple groups.
Mainly, omit Part~C and  replace Part~E by the following:

Let $M$ be the normal subgroup of $G$ such that $M/A = K$.
By Lemma~\ref{sylow},
$C$ normalizes some nontrivial Sylow subgroup $U$ of $M$.
By Frattini argument, $U$ is intravariant in $G$.
\end{remark}

\section{Profinite Groups}

\begin{theorem}\label{main}
Let $G$ be a $C$-group 
and let $G \to F$ be a $C$-epimorphism 
onto a $C$-Frobenius group $F$.
Then $G$ contains
a $C$-Frobenius subgroup $H$
that maps onto $F$.
Moreover,
we may assume that $F$ and $H$ have precisely the same prime divisors.
\end{theorem}

\begin{proof}
The last assertion of the theorem follows immediately from the rest,
because we may drop from the Frobenius kernel of $H$
its $p$-primary components for those primes $p$
that do not divide the order of the Frobenius kernel of $F$.

If necessary,
replace $G$ by a closed subgroup 
to assume that $G$ is finitely generated.
Then there is a sequence 
$N_0 \ge N_1 \ge N_2 \ge \ldots$ 
of open normal subgroups of $G$ such that
$\bigcap_i N_i = \{1\}$.
Let $A = \ker(G \to F)$.
Since $\bigcap_i N_i A = A$,
without loss of generality
$C \cap N_i A = \{1\}$ for every $i$.

Put $M_i = N_i \cap A$ for every $i$.
Then $F = G/M_0$ and $G = \inverselimit_i G/M_i$.
For each $i$ we have the following cartesian diagram
of epimorphisms of $C$-groups
$$
\xymatrix{%
G/M_{i+1} \ar[r] \ar[d] & G/N_{i+1} \ar[d] \\
G/M_{i} \ar[r] & G/M_i N_{i+1} \\
}%
$$
in which the groups on the right handed side are finite.
By induction hypothesis
$G/M_i$ contains 
a $C$-Frobenius subgroup $F_i$
that maps onto $F$.
Its image in $G/M_iN_{i+1}$ is a $C$-Frobenius subgroup.
By Theorem~\ref{main finite}
it lifts to a $C$-Frobenius subgroup of $G/N_{i+1}$.
Hence by Lemma~\ref{fiber product}
$F_i$ lifts to a $C$-Frobenius subgroup $F_{i+1}$ of $G/N_{i+1}$.
Thus 
$H = \inverselimit_i F_i$ is 
a $C$-Frobenius subgroup of $G$ that maps onto $F$.
\end{proof}

\begin{corollary}
Let $\beta \colon B \to C$ be an epimorphism
of a finite group $B$ onto a finite nontrivial cyclic group $C$
and let $F$ be a $C$-Frobenius group
with Frobenius kernel
$K \cong \hat{\mathbb Z}_\pi$
for some set $\pi$ of primes.
Then there is a $C$-embedding
$F \to B \coprod C$.
\end{corollary}

\begin{proof}
By Theorem~\ref{main}
it suffices to construct a $C$-epimorphism
$B \coprod C \to F$.
Let $k$ be a generator of $K$.
The epimorphism $B \to C$ given by
$b \mapsto \beta(b)^k$
together with the identity map $C \to C$
define a $C$-homomorphism $B \coprod C \to F$.
Its image $\langle C, C^k \rangle$
contains $C$ and 
$(c^k)^{-1} c = k^{-1} k^c \in K$,
where $c$ is a generator of $C$.
Since $F = CK$,
it suffices to show that
$k^{-1} k^c$ generates $K$.

If this assertion holds for each Sylow subgroup of $K$
(and its generator instead of $k$)
then it holds for $K$.
Thus we may assume that $K \cong {\mathbb Z}_p$
for some prime $p$.
Finally, we may replace $K$ by its Frattini quotient.
$\bar K \cong {\mathbb Z}/p{\mathbb Z}$.

As $c \ne 1$,
we have 
$k^{-1} k^c \ne 1$ and hence $k^{-1} k^c$ generates $K$.
\end{proof}

\begin{example}\label{Pop's question}
Let $K = {\mathbb Z}_7$.
Then $C = {\mathbb Z}/6{\mathbb Z}$
acts fixed-point freely on $K$
and hence the semidirect product $F = CK$
is a $C$-Frobenius group.
By the corollary
there is an embedding $\lambda \colon F \to G_1 \coprod G_2$,
where $G_1 = G_2 = C$,
such that $\lambda(C) = G_2$.
Let $H = \lambda(F)$.
Then 
$H$ is prosolvable, but
\begin{enumerate}
\item
$H \cap G_2 = G_2$ is of order $6$,
and hence is not an $\ell$-group for some prime $\ell$;
\item
$H$ is infinite,
and hence
$H \le G_i^\sigma$ for no $i \in I$ and no $\sigma \in G$.
\end{enumerate}
This answers the question of Pop mentioned in the introduction.
\end{example}

\end{document}